\documentclass{conm-p-l}

\newtheorem{theorem}{Theorem}[section]
\newtheorem{lemma}[theorem]{Lemma}

\theoremstyle{definition}
\newtheorem{definition}[theorem]{Definition}

\theoremstyle{remark}
\newtheorem{remark}[theorem]{Remark}

\numberwithin{equation}{section}

\usepackage{amssymb}
\usepackage{amsmath}
\usepackage{latexsym}
\usepackage[mathscr]{euscript}
\usepackage{graphicx}
\usepackage{amscd}
\usepackage{amsthm} 
\usepackage{enumerate}
\usepackage{wrapfig}
\newtheorem{proposition}[theorem]{Proposition}
\newtheorem{conjecture}[theorem]{Conjecture}
 
 \newtheorem{question}[theorem]{Question}

\newcommand{\bz}{\mathbb{Z}}
\newcommand{\br}{\mathbb{R}}

\newcommand{\bh}{\mathbb{H}}

\newcommand{\p}{\partial}
\newcommand{\cc}{\mathcal{C}}

\newcommand{\cs}{\mathcal{S}}

\newcommand{\cp}{\mathcal{P}}

\newcommand{\cm}{\mathcal{M}}

\newcommand{\cf}{\mathcal{F}}

\newcommand{\ca}{\mathcal{A}}

\newcommand{\hk}{\hookrightarrow}

\newcommand{\la}{\longrightarrow}
\newcommand{\bfl}{\begin{flushleft}}
\newcommand{\efl}{\end{flushleft}}

\newcommand{\sma}{\wedge}

 \newcommand{\co}{\mathcal{O}}
 \newcommand{\cd}{\mathcal{D}}

\newcommand{\xr}{\xrightarrow}

\begin{document}

 \title{Open-closed field theories,  string topology, and Hochschild homology}

 \author {Andrew J. Blumberg}
 \address{Dept. of Mathematics,  Stanford University,  
 Stanford, CA 94305}
 \email{blumberg@math.stanford.edu}
  \thanks{The first author was partially supported by an NSF postdoctoral fellowship}
 
\author{Ralph L. Cohen}
\address{Dept. of Mathematics,  Stanford University,  
 Stanford, CA 94305}
 \email{ralph@math.stanford.edu}
  \thanks{The second author was  partially supported by   NSF grant  DMS-0603713}

\author{Constantin Teleman}
\address{Dept. of Mathematics, 
 UC Berkeley, Berkeley, CA 94720}
 \email{teleman@math.berkeley.edu}
 \thanks{The third author was  partially supported by a research grant from the NSF} 
 
 \date{\today}
  
\begin{abstract} 
  In this expository paper we discuss a project regarding the string
  topology of a manifold, that was inspired by recent work of
  Moore-Segal \cite{mooresegal}, Costello \cite{costello}, and Hopkins
  and Lurie \cite{lurie} on ``open-closed topological conformal field
  theories".  Given a closed, oriented manifold $M$, we
  describe the ``string topology category" $\cs_M$, which is enriched
  over chain complexes over a fixed field $k$.  The objects of $\cs_M$
  are connected, closed, oriented submanifolds $N$ of $M$, and the
  complex of morphisms between $N_1$ and $N_2$ is a chain complex
  homotopy equivalent to the singular chains $C_*(\cp _{N_1, N_2})$
  where $\cp_{N_1, N_2}$, is the space of paths in $M$ that start in
  $N_1$ and end in $N_2$.  The composition pairing in this category is
  a chain model for the open string topology operations of Sullivan
  \cite{sullivan}, and expanded upon by Harrelson \cite{harrelson} and
  Ramirez \cite{ramirez}.  We will describe a calculation yielding
  that the Hochschild homology of the category $\cs_M$ is the homology
  of the free loop space, $LM$.  Another part of the project is to
  calculate the Hochschild cohomology of the open string topology
  chain algebras $C_*(\cp_{N,N})$ when $M$ is simply connected, and
  relate the resulting calculation to $H_*(LM)$.  These calculations
  generalize known results for the extreme cases of $N= point$ and $N
  = M$, in which case the resulting Hochschild cohomologies are both
  isomorphic to $H_*(LM)$.  We also discuss a spectrum level analogue
  of the above results and calculations, as well as their relations to
  various Fukaya categories of the cotangent bundle $T^*M$ with its
  canonical symplectic structure.  This paper is purely expository,
  and is intended to be a background survey and announcement of some
  of our recent results. Details will appear in due course
  \cite{bct}.  
\end{abstract}

\maketitle

 \tableofcontents

 \section*{Introduction}  In an open-closed topological field theory, one studies cobordisms between compact one-dimensional manifolds, whose boundary components are labeled by an indexing set, $\cd$. The cobordisms are those of manifolds with boundary, that preserve the labeling sets in a specific way.  The set of labels $\cd$  are referred to as  ``D-branes",  and  in the string theory literature these are boundary values of ``open strings".    An open-closed field theory is a monoidal functor
 from a category built out of such manifolds and cobordisms, that takes values in  a linear category, such as vector spaces,  chain complexes, or even  the  category of spectra.
 In this paper we will discuss two flavors of such open-closed field theories: ``topological quantum field theories" (TQFT)  as introduced by Moore and Segal \cite{mooresegal}, and ``topological conformal field theories", (TCFT),  as studied by Getzler \cite{getzler} and Costello \cite{costello}. 
 
 The open part of such a theory $\cf$  is the restriction of  $\cf$ to the ``open subcategory".   This  is the full subcategory generated by  those compact one-manifolds, all of whose path components have nonempty boundary.  As Moore and Segal originally pointed out,  the data of an open field theory can be encoded in a category (or as Costello points out, an $A_\infty$-category when $\cf$ is an open-closed TCFT),  $\cc_\cf$.  The objects of $\cc_\cf$ are the set of $D$-branes, $\cd$. The space of morphisms
 between $\lambda_0$ and $\lambda_1 \in \cd$ is given by the value of the theory $\cf$ on the object $I_{\lambda_0, \lambda_1}$, defined by the interval $[0,1]$ where the boundary component
 $0$ is labeled by $\lambda_0$, and $1$ is labeled by $\lambda_1$.  We denote this vector space by $\cf (\lambda_0, \lambda_1)$.   The composition rules in this ($A_\infty$) category are defined by the values of $\cf$ on certain ``open-closed" cobordisms.  Details of this construction will be given below.
 
 In this paper we will report on a project whose goal is to understand how the ``String Topology" theory of a manifold fits into this structure.  This theory, as originally introduced by Chas and Sullivan \cite{chassullivan} starts with a closed, oriented $n$-dimensional manifold $M$.  It was shown in \cite{cohengodin}
 that there is a (positive boundary) TQFT $\cs_M$, which assigns to a circle the homology of the free loop space,
 $$
 \cs_M(S^1) = H_*(LM; k)
 $$
 with field coefficients.  This was recently extended by Godin \cite{godin} to show that string topology is actually an open-closed homological conformal field theory.  In this theory the set of D-branes $\cd_M$ is the set of connected, closed, oriented, connected submanifolds of $M$.  The theory assigns to a compact one-manifold
$c$ with boundary levels, the homology of the mapping space,
$$
\cs_M(c) = H_*(Map (c, \p;  M)).
$$
Here $Map (c, \p; M)$ refers to the space of maps $c \to M$ that take the labeled boundary components to the submanifolds determined by the labeling.
In particular, we write $\cp_{N_0, N_1} = Map (I_{N_0, N_1}, \p; M)$ for  the space of paths   $\gamma: [0,1]  \to M$ such that $\gamma (0) \in N_0$, and $\gamma (1) \in N_1$.   In Godin's theory, given any open-closed cobordism $\Sigma_{c_1, c_2}$ between one-manifolds $c_1$ and $c_2$, there are homological operations
$$
\mu_{\Sigma_{c_1, c_2}} : H_*(BDiff (\Sigma_{c_1, c_2}); k) \otimes H_*(Map (c_1, \p; M); k) \la H_*(Map (c_2, \p; M); k).
$$

An open-closed topological conformal field theory in the sense of
Costello is a \sl chain complex \rm valued theory, and it is
conjectured that the string topology theory has the structure of such
a theory.  The following theorem, which we report on in this paper,
gives evidence for this conjecture.

\begin{theorem}\label{bct1}
Let $k$ be a field.
\begin{enumerate}
\item {There exists a $DG$-category (over $k$) $\cs_M$, with the following properties:
\begin{enumerate}
\item The objects are the set of $D$-branes, $\cd_M$   =  \{connected, closed, oriented submanifolds of $M$\} 
\item The morphism  complex, $Mor_{\cs_M}(N_1, N_2)$, is  chain homotopy equivalent to the singular chains on the path space $C_*(\cp_{N_1, N_2})$.
\end{enumerate}
The compositions in $\cs_M$ realize the open-closed string topology operations on the level of homology.
}

\item The Hochschild homology of $\cs_M$ is equivalent to the homology of the free loop space, 
$$
HH_*(\cs_M)  \cong H_*(LM; k).
$$
\end{enumerate}

\bf Note.  \rm In this theorem we construct a $DG$-category with
strict compositions rather than an $A_\infty$ category.  See section 2
below.
\end{theorem}

Given any fixed submanifold $N$, the space of self-morphisms,
$Mor_{\cs_M}(N, N) \simeq C_*(\cp_{N,N})$ is a differential graded
algebra. Again, on the level of homology, this algebra structure is
the string topology product introduced by Sullivan \cite{sullivan}.
In this note we pose the following question and report on its answer
in a variety of special cases.  (See Theorem \ref{doublecent} below.)
Details will appear in \cite{bct}.
 
\begin{question}\label{onebrane}
Let  $M$ be a simply connected, closed submanifold.  For which connected,  oriented, closed  submanifolds   $N \subset M$  is the  Hochschild cohomology of $C_*(\cp_{N,N})$ isomorphic to the  homology of the free loop space,
$$
HH^*(C_*(\cp_{N,N}),  C_*(\cp_{N,N})) \cong H_*(LM) 
$$
as algebras?  The algebra structure of the left hand side is given by cup product in Hochschild cohomology, and on the right hand side by the  Chas-Sullivan   string topology product. 
\end{question}
 
We observe that in the two extreme cases ($N$ a point, and $N = M$),
affirmative answers to this question are known.  For example, when $N$
is a point, $\cp_{N,N}$ is the based loop space, $\Omega M$, and the
statement that $HH_*(C_*(\Omega M), C_*(\Omega M)) \cong H_*(LM)$ was
known in the 1980's by work of Burghelea, Goodwillie, and others.  The
Hochschild cohomology statement then follows from Poincar\'{e} duality.
Similarly, when $N = M$, then $ \cp_{N,N} \simeq M$, and the string
topology algebra on $C_*(\cp_{N,N})$ corresponds, via Poincar\'{e}
duality, to the cup product in $C^*(M)$.  The fact that the Hochschild
cohomology of $C^*(M)$ is isomorphic to $H_*(LM)$ follows from work of
J. Jones in the 1980's, and the fact that the ring structure
corresponds to the Chas-Sullivan product was proved in
\cite{cohenjones}.  In this note we are able to report on a
calculation of $HH^*(C_*(\cp_{N,N}), C_*(\cp_{N,N}))$ which yields an
affirmative answer to this question in many cases. (See Theorem
\ref{doublecent} below.)  These cases include when the inclusion map
$N \hk M$ is null homotopic.
Thus
$$ HH^*(C_*(\cp_{N,N}),  C_*(\cp_{N,N})) \cong H_*(LS^n)$$
for every connected, oriented, closed submanifold of a sphere $S^n$.
Other cases when one gets an affirmative answer to the above question
include when the inclusion $N \hk M$ is the inclusion of the fiber of
a  fibration $p : M \to B$, or more generally, when $N \hk M$ can be
factored as a sequence of embeddings, $N =N_0 \hk N_1 \hk \cdots N_{i}
\hk N_{i+1} \cdots N_k = M$ where each $N_i \subset N_{i+1}$ is the
inclusion of a fiber of a fibration $p_{i+1} : N_{i+1} \to B_{i+1}$.        

    We point out that an amusing aspect of this question is that for
    any $N \hk M$ for which the answer is affirmative, then one can
    use this submanifold as a single $D$-brane and recover $H_*(LM)$
    as a Hochschild cohomology ring (i.e., ``one brane is enough"),
    and that all such branes yield the same answer.
 
 This paper is organized as follows. In section one below we discuss
 the relevant background of open-closed topological field theories,
 including a review of work of Moore and Segal \cite{mooresegal}, and
 of Costello \cite{costello}.  In section 2 we describe the
 ingredients of the proof of Theorem~\ref{bct1} and discuss the
 Hochschild cohomology calculations of the chain algebras,
 $C_*(\cp_{N,N})$ in Theorem \ref{doublecent} below.  The methods
 involve generalized Morita theory, and so yield comparisons between
 certain module categories over the algebras $C_*(\cp_{N,N})$.  We present these in Theorem \ref{module} below.  In the 
 extreme cases mentioned above, these comparisons reduce to the
 standard equivalences of certain module categories over the cochains
 $C^*(M)$ and the chains of the based loop space, $C_*(\Omega M)$
 (originally obtained in \cite{DGI}). 
 In section 3 we discuss possible relationships between the categories
 described here and certain Fukaya categories of the cotangent bundle,
 $T^*M$ with its canonical symplectic structure.

{\bf Acknowledgments}

We would like to thank Bill Dwyer,  Michael
Hopkins, and  Michael Mandell   for their help with various aspects of this project.

\section{Open-closed Topological  Field Theories}     As mentioned in the introduction,    the objects  of study in an open-closed field theory  are parameterized, compact, oriented one-manifolds, $c$, together with a labeling of the  components of the boundary, $\p c$, by elements of a set, $\cd$.    An  ``open-closed" cobordism $\Sigma_{c_1, c_2}$ between two objects $c_1$ and $c_2$  is an oriented surface $\Sigma$, whose boundary is partitioned into three parts:  the incoming boundary, $\p_{in} \Sigma$ which is identified with $c_1$,  the outgoing boundary $\p_{out}\Sigma$ which is identified with $c_2$, and the remaining part of the boundary, referred to as the ``free part",  $\p_{free} \Sigma$ whose path components are labeled by $\cd$,  with the property that
$\p_{free} \Sigma$ is itself a cobordism between $\p c_1$ and $\p c_2$,  preserving the labeling.  This is the usual notion of a cobordism of manifolds with boundary,  with the additional data of the labeling set $\cd$.  Figure 1 below  is a picture of a one-manifold whose boundary components are labeled by elements of $\cd$, and figure 2 is a picture of an open-closed cobordism.  In this picture the free part of the boundary, $\p_{free} \Sigma$ is highlighted in red.  In figure 3 a smooth surface is shown that is homeomorphic to the open-closed cobordism
given in figure 2.  The free part of the boundary is again highlighted in red.

\begin{figure}[ht]
  \centering
  \includegraphics[height=8cm]{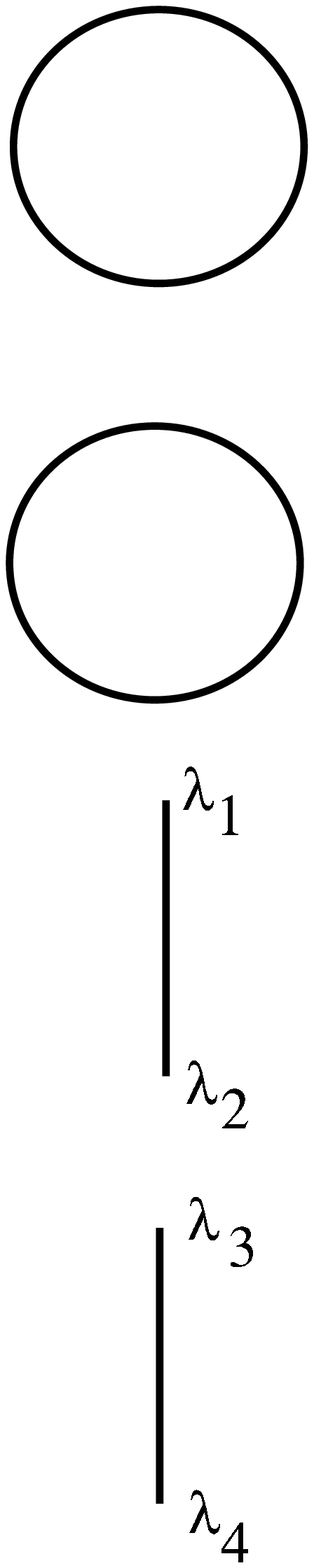}
 \caption{ A one manifold with labels $\lambda_i \in \cd$}   
 \label{figone}
\end{figure}

\begin{figure}[ht]
  \centering
  \includegraphics[height=4cm]{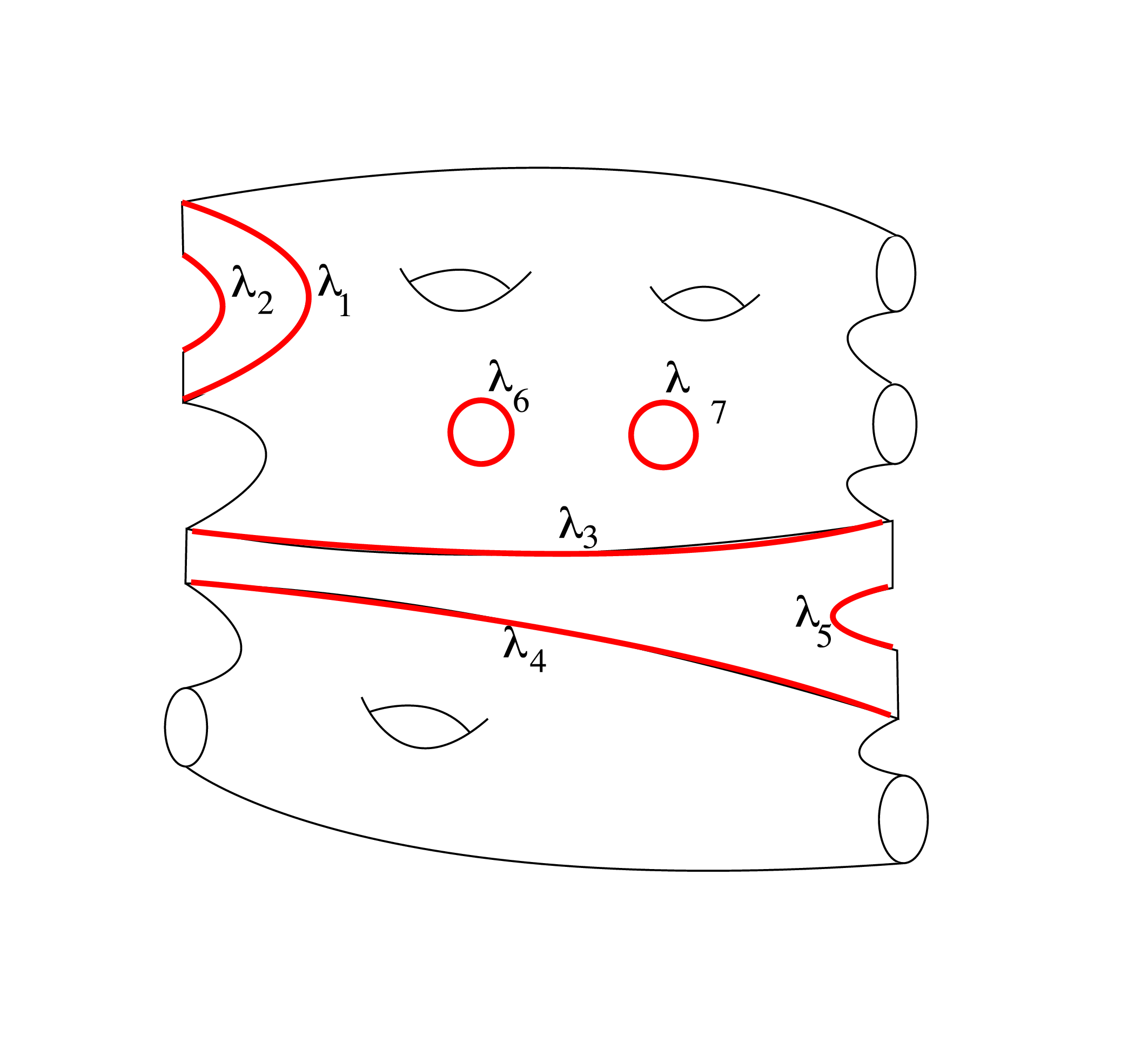}
 \caption{An open-closed cobordism }
 \label{figtwo}
\end{figure}
 
\begin{figure}[ht]
  \centering
  \includegraphics[height=4cm]{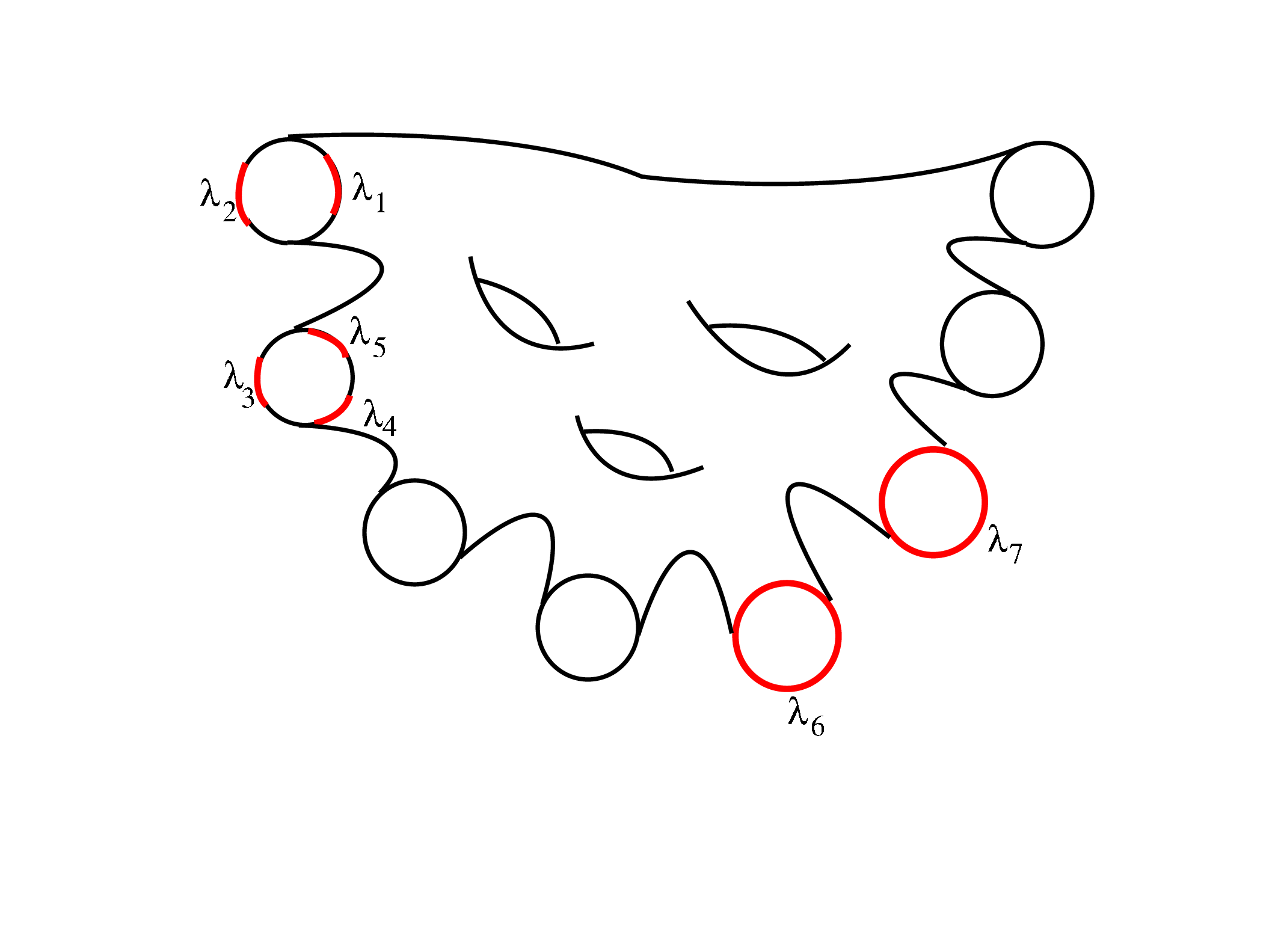}
  \caption{A smooth  open-closed cobordism }
  \label{figthree}
\end{figure}

In \cite{mooresegal},  Moore and Segal describe basic properties of \sl open-closed topological quantum field theories,  \rm  and in a sense,
Costello then gave a derived version of this theory when he gave a description of \sl open-closed topological conformal field theories. \rm

\subsection{The work of Moore and Segal on open-closed TQFT's}

In \cite{mooresegal} Moore and Segal describe how   an \sl open-closed topological quantum field theory \rm $\cf$  assigns to each one manifold $c$ with boundary components labeled by $\cd$, a vector space over a field $k$,  $\cf (c)$.   The theory $\cf$ also  assigns to every diffeomorphism class of open-closed cobordism $\Sigma_{c_1, c_2}$ a linear map
$$
\cf (\Sigma_{c_1, c_2} ) : \cf (c_1) \to \cf (c_2).
$$

This assignment is required to satisfy two main properties:

\begin{enumerate}
\item 1.  \sl Gluing: \rm  One can glue two open-closed cobordisms when the outgoing boundary of one is identified with the incoming boundary of the other:
$$
\Sigma_{c_1, c_2} \# \Sigma_{c_2, c_3} = \Sigma_{c_1, c_3}.
$$
In this case the operation $\cf (\Sigma_{c_1, c_2} \# \Sigma_{c_2, c_3})$ is required to be the composition:
$$
\cf (\Sigma_{c_1, c_2} \# \Sigma_{c_2, c_3}) : \cf(c_1) \xr{\cf (\Sigma_{c_1, c_2})} \cf (c_2)  \xr{\cf (\Sigma_{c_2, c_3})} \cf (c_3).
$$
This condition can be viewed as saying that $\cf$ is a functor $\cf :\cc_\cd  \to Vect_k$, where    $\cc_\cd$ is the cobordism category whose objects are one manifolds with boundary labels in $\cd$, and whose morphisms are diffeomorphism classes of open-closed cobordisms.   Here the diffeomorphisms are required to preserve the orientations, as well as the boundary structure ($\p_{in},  \p_{out}$, and the labeling).  $Vect_k$ is the category of vectors spaces over the field $k$,  whose morphisms are linear transformations between them. 
 \item \sl Monoidal:  \rm There are required to be natural isomorphisms,
 $$
 \cf(c_1) \otimes \cf(c_2) \xr{\cong} \cf(c_1 \sqcup c_2)
 $$
 that makes $\cf$ into a monoidal functor.  (The monoid structure in $\cc_\cd$ is given by disjoint union of both the object manifolds and the morphism cobordisms. )
 \end{enumerate}

Let $\cf$ be an open-closed TQFT.  The value of $\cf$ on a single unit circle, $\cf (S^1)$ is known as the \sl closed state space \rm of the theory $\cf$.  It is well known that $\cf$ is a  commutative \sl Frobenius algebra \rm over $k$.   That is, there is an associative multiplication $\mu_{\cf} : \cf (S^1) \otimes \cf (S^1) \to \cf (S^1)$   coming from the value of $\cf$ on the pair of pants cobordism from $S^1 \sqcup S^1$ to $S^1$.  The unit disk, viewed as having one  outgoing boundary component, $\p_{out} D^2 = S^1$, is a cobordism
from the emptyset $\emptyset$  to $S^1$,  and   therefore induces a map $\iota: k \to \cf (S^1)$, which is the unit in the algebra structure.  Thinking of the disk as having one \sl incoming  \rm boundary component, $\p_{in} D^2 = S^1$,   induces a map $\theta_{\cf} : \cf (S^1) \to k$ which is the ``trace map" in the theory.  That is,
the bilinear form
$$\langle \, , \, \rangle : \cf(S^1) \times \cf (S^1)   \xr{\mu_{\cf}} \cf(S^1) \xr{\theta_{\cf}} k$$ is nondegenerate.

\begin{figure}[ht]
  \centering
  \includegraphics[height=4
  cm]{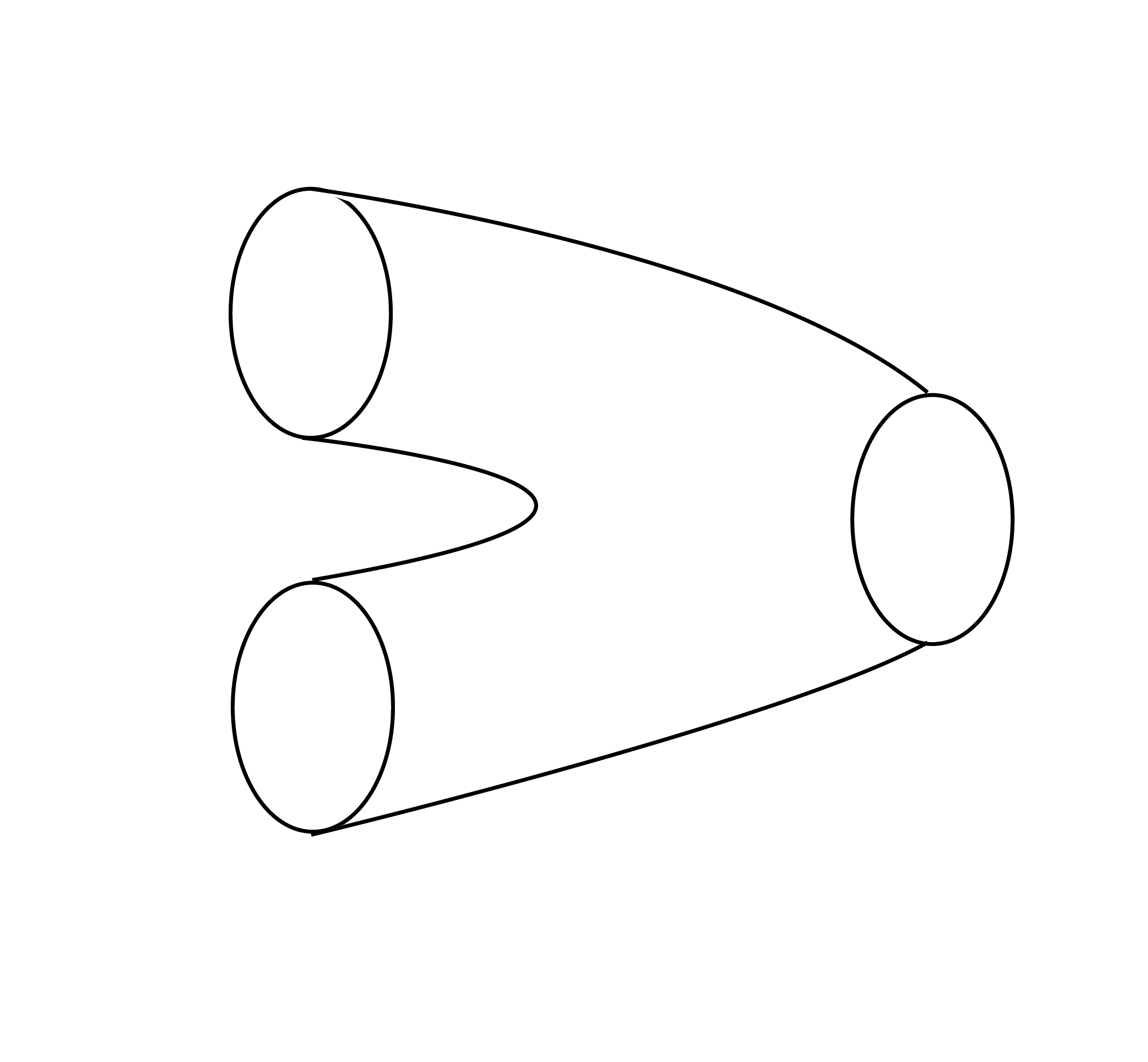}
  \caption{The pair of pants cobordism inducing the multiplication $\mu_{\cf}: \cf(S^1) \otimes \cf (S^1) \to \cf (S^1)$. }
  \label{figfour}
\end{figure}

There is more algebraic structure associated to an open-closed field theory $\cf$.  As described by Moore and Segal,    there is a category,
$\cc_\cf$ associated to the open part of the field theory.

\begin{definition}\label{category}   The category $\cc_{\cf}$ associated to an open-closed TQFT $\cf$ has as its objects
  the set of $D$-branes, $\cd$.  The space of morphisms  between objects $\lambda_1$ and $\lambda_2$ is given by the value of  the field theory $\cf$ on the one-manifold $I_{\lambda_1, \lambda_2}$ which is given by the interval  $[0,1]$  with boundary components  labeled by $\lambda_1$ and $\lambda_2$.  We write this space as $\cf(\lambda_1, \lambda_2)$.  

\begin{figure}[!]
  \centering
  \includegraphics[height=6cm]{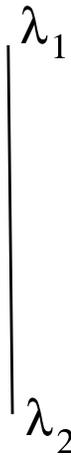}
  \caption{The  one manifold $I_{\lambda_1, \lambda_2}$ which induces the morphism space $\cf(\lambda_1, \lambda_2)$. }
  \label{figfive}
\end{figure}

The composition law in the category $\cc_\cf$ is defined by the open-closed cobordism shown in figure 6. 

\begin{figure}[ht]
  \centering
  \includegraphics[height=5cm]{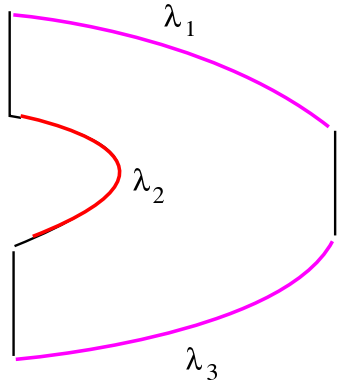}
  \caption{The value of $\cf$ on this cobordism defines the composition pairing, $\cf (\lambda_1, \lambda_2) \otimes \cf (\lambda_2, \lambda_3 )\to \cf (\lambda_1, \lambda_3)$}
  \label{figsix}
\end{figure}

\end{definition}

Notice that the endomorphism algebras in this category, $\cf (\lambda, \lambda)$, are also Frobenius algebras. For simplicity we write these algebras as $\cf (\lambda)$.   The trace maps are induced by the open-closed cobordism between $I_{\lambda, \lambda}$ and the empty set given by the disk as in figure 7.

\begin{figure}[ht]
  \centering
  \includegraphics[height=4cm]{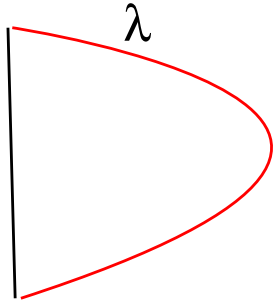}
  \caption{The value of $\cf$ on this open-closed cobordism defines the  trace map $\cf (\lambda)  \to k$}
  \label{figseven}
\end{figure}

We observe that the closed state space $\cf (S^1)$ is necessarily commutative as an algebra, because  the cobordisms shown in figure 8  admit an orientation preserving
diffeomorphism between them  that fixes the boundary pointwise.      However for a $D$-brane $\lambda \in  \cd$,  the fact that the open-closed cobordisms shown in figure 9  \sl are not \rm diffeomorphic via an orientation preserving diffeomorphism that fixes the incoming and outgoing boundaries,  imply that the Frobenius algebra $\cf (\lambda, \lambda)$ may not be  commutative.

\begin{figure}[ht]
  \centering
  \includegraphics[height=4cm]{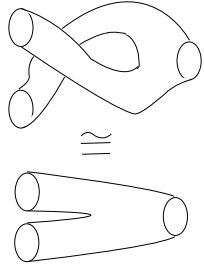}
  \caption{These diffeomorphic cobordisms imply that the Frobenius algebra $\cf (S^1)$ is commutative.}
  \label{figeight}
\end{figure}

\begin{figure}[ht]
  \centering
  \includegraphics[height=4cm]{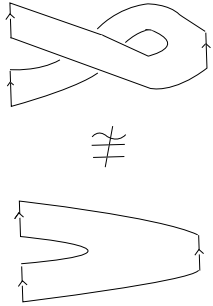}
  \caption{These surfaces are not diffeomorphic as open-closed cobordisms, and thus the Frobenius algebras $\cf (\lambda)$ may not be commutative.}
  \label{fignine}
\end{figure}

These   algebras are, of course, related to each other.  For example, the ``whistle"  open-closed cobordism given in figure 10 defines a ring homomorphism
$\theta_\lambda : \cf (S^1)   \to  \cf (\lambda)$, which, is easy to see takes values in the center $Z(\cf (\lambda))$ (see \cite{mooresegal} for details.)
So in particular one has the following result.

\begin{proposition}\label{center}  Any open-closed TQFT $\cf$  comes equipped with  map of algebras
$$\theta_\lambda : \cf(S^1) \to Z(\cf (\lambda))$$
where $Z(\cf (\lambda))$ is the center of the endomorphism algebra $\cf (\lambda)$,  for any $\lambda \in \cd$.     
\end{proposition}

\begin{figure}[ht]
  \centering
  \includegraphics[height=3 cm]{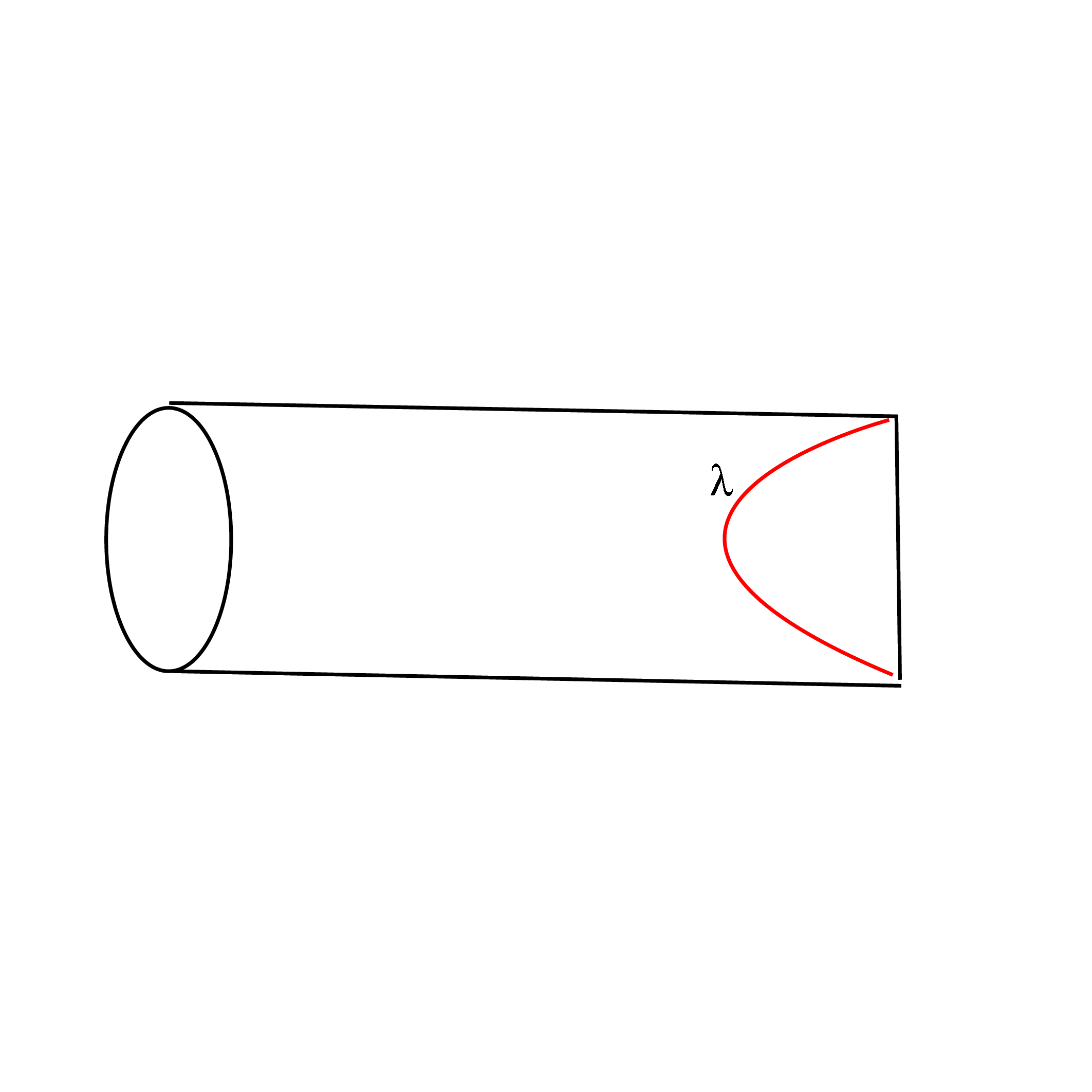}
  \caption{The ``whistle open-closed cobordism"  inducing the map $\theta_\lambda: \cf(S^1) \to Z(\cf (\lambda)).$  }
  \label{figten}
\end{figure}

Turning the whistle cobordism around,  so that its incoming boundary is $I_{\lambda, \lambda}$, and its outgoing boundary is $S^1$,
defines a homomorphism $\theta^*_\lambda : \cf (S^1) \to \cf (\lambda)$, which is not difficult to see is adjoint to $\theta_\lambda$, with respect to the inner products defined by the corresponding Frobenius algebras.  Moreover, studying the relevant glued cobordisms, one can show that the composition,
$\theta_\lambda \circ \theta^*_\lambda$ satisfies the ``Cardy formula",
\begin{equation}\label{cardy}
\theta_\lambda \circ \theta^*_\lambda (\phi) = \sum_{i = 1}^n \psi^i \phi \psi_i
\end{equation}
where $\{\psi_1, \cdots , \psi_n\}$ is any basis of $\cf (\lambda)$, and $\{\psi^1, \cdots , \psi^n\}$ is the dual basis (with respect to the inner product in the Frobenius algebra structure).   Again, see \cite{mooresegal} for the details of this claim.

\subsection{The work of  Costello on open-closed TCFT's}
 
\newcommand{\ocd}{\co \cc_{\cd}}
\newcommand{\od}{\co_{\cd}}
In \cite{costello} Costello studied open-closed \sl topological conformal field theories (TCFT).  \rm  Such a theory can be viewed as a derived version of a topological quantum field theory,
and in a sense, Costello's work can,  in part be viewed as a derived extension and generalization of the work of Moore and Segal.

More precisely,  the TCFT's Costello studied are functors,    
$$
\cf : \ocd \to Comp_k
$$  where $\ocd$ is an open-closed cobordism category, enriched over chain complexes,  and $Comp_k$  is the symmetric monoidal category of chain-complexes over a ground field $k$.   In Costello's work, $Char (k) = 0$.     By an ``open-closed cobordism category, enriched over chain complexes", Costello means the following.  Let $\cd$ be  an indexing set  of  ``D-branes" as above.  Then the objects of $\ocd$ are parameterized, compact, oriented one-manifolds, $c$, together with a labeling of the  components of the boundary, $\p c$, by elements of $\cd$, as described in the previous section.

  To describe  the chain complex  of  morphisms between objects  $c_1$ and $c_2$,  one considers the   moduli space of all
Riemann surfaces  that form open-closed cobordisms between $c_1$ and $c_2$.   This moduli space was originally described by Segal \cite{segalconf}
when the $c_i$'s have no boundary.  For the general situation we refer the reader to Costello's paper \cite{costello}.
These open-closed cobordisms are required to satisfy the additional
``positive boundary" requirement, that every path component of an element $\Sigma \in \cm_{\cd} (c_1, c_2)$ has a nonempty incoming
boundary.    It is standard to see that
$$
\cm_{\cd} (c_1, c_2) \simeq  \coprod BDiff^+ (\Sigma, \p \Sigma)
$$
where the disjoint union is taken over all diffeomorphism classes of open-closed cobordisms from $c_1$ to $c_2$.  These diffeomorphisms are diffeomorphisms
of open-closed cobordisms, as defined in the previous section.  They 
are orientation preserving,  they preserve the  incoming and outgoing boundaries pointwise, and  they preserve  the labelings in 
$\cd$.    The morphisms in $\ocd$ are then the singular chains with coefficients in $k$, $Mor_{\ocd}(c_1, c_2) = C_*(\cm_{\cd} (c_1, c_2); k)$.

A topological conformal field theory is then a functor $\cf : \ocd \to Comp_k$ which is ``$h$-monoidal", in the sense that there
are natural transformations
$\cf(c_1) \otimes \cf(c_2) \to \cf (c_1 \sqcup c_2)$ which are quasi-isomorphisms of chain complexes. Costello calls $\ocd - mod$ the functor category of topological conformal field theories.

Let  $\od \hk \ocd$ be the full subcategory whose objects have no closed components.  That is, every connected component of a one-manifold $c \in Ob (\od)$ has (labeled)  boundary.  Write $\od-mod$ to be the functor category of $h$-monoidal functors $\phi : \od \to Comp_k$.  We refer to such a functor as an ``open-field theory".  

Costello observed that an open topological conformal field theory $\phi : \od \to Comp_k$ defines an $A_\infty$-category, enriched over chain complexes,  in much the same was as an open topological  quantum field theory defines a category (see Definition (\ref{category}) above). This is most easily seen if the field theory
is \sl strictly \rm monoidal, that is, the transformations $ \phi (c_1) \otimes \phi (c_2) \to \phi (c_1 \sqcup c_2)$ are \sl isomorphisms \rm of chain complexes, rather than only quasi-isomorphisms.  In this case the associated $A_\infty$-category, which we call $ \cc_\phi $, has objects given
by the set of $D$-branes $\cd$.  The  space of   morphisms $\phi (\lambda_0, \lambda_1)$ is the chain complex given by the value of the functor
$\phi$ on the object $I_{\lambda_0, \lambda_1}$.  We call this space $\phi (\lambda_0, \lambda_1)$.     
   The higher  compositions 
$$
\phi (\lambda_1, \lambda_2) \otimes \phi (\lambda_2, \lambda_2) \otimes \cdots \otimes \phi (\lambda_{n-1}, \lambda_n) \la \phi (\lambda_1, \lambda_n)
$$
 are given by  the value of the functor $\phi$ on the open-closed cobordism  between $\coprod_{i=1}^{n-1} I_{\lambda_i, \lambda_{i+1}}    $ and  $I_{\lambda_1, \lambda_n}$ given by the connected, genus zero surface $D_{\lambda_1, \cdots \lambda_n}$  pictured  in figure 11 in the case $n = 4$.  

\begin{figure}[ht]
  \centering
  \includegraphics[height=6cm]{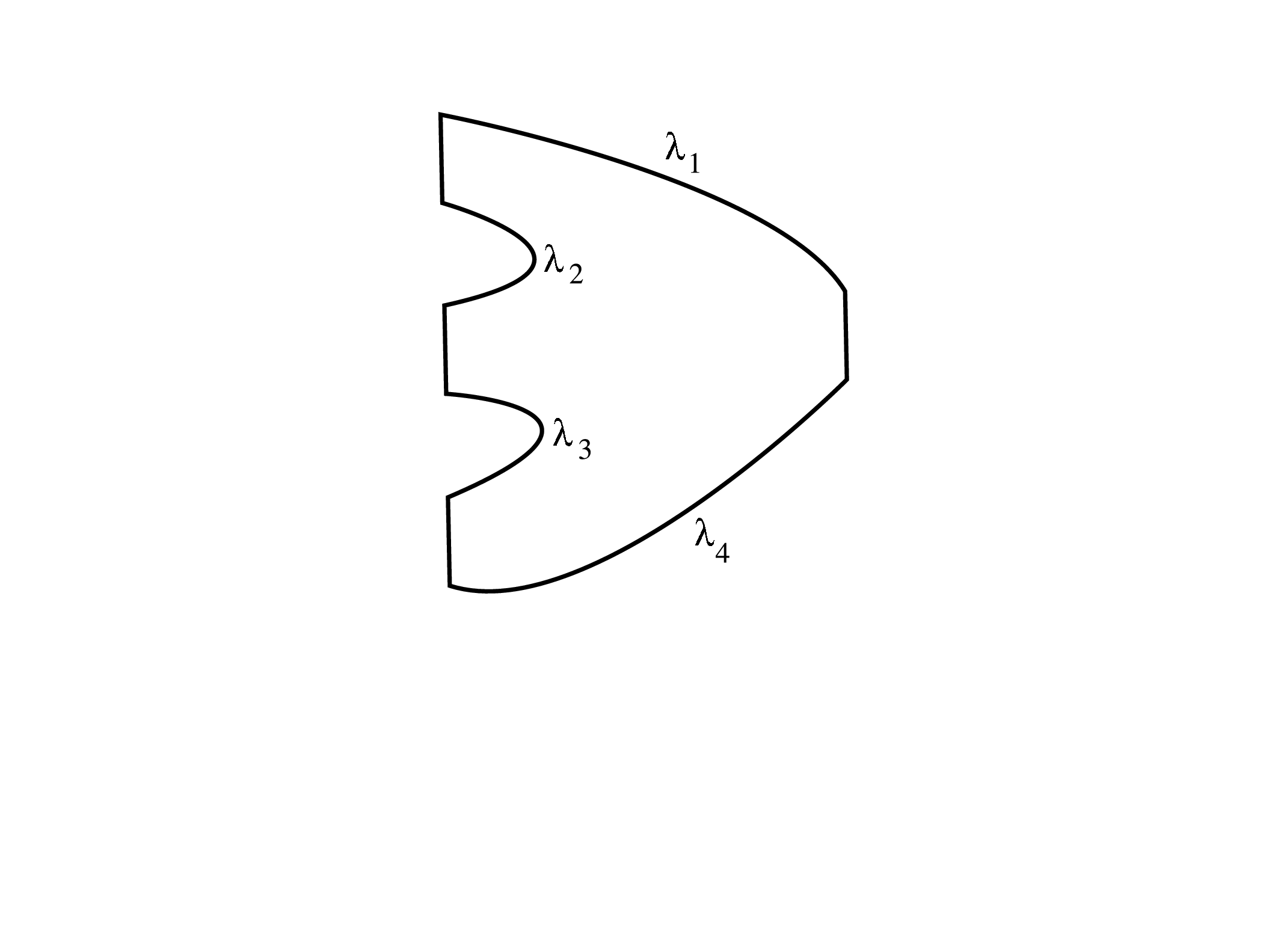}
 \caption{The open-closed cobordism $D_{\lambda_1, \cdots \lambda_4}$ }
 \label{figeleven}
\end{figure}

 The $A_\infty$-category defined by an open TCFT  has additional properties that Costello referred to as a ``Calabi-Yau" $A_\infty$ category.  The following theorem of Costello describes the central nature of this category in open-closed field theory.
 
 \begin{theorem}\label{costello}(Costello \cite{costello}) a.   The 
  restriction functor $\rho : \ocd-mod \to \od-mod$  from open-closed TCFT's to open TCFT's has a derived left adjoint,
 $L_\rho  : \od-mod \to \ocd- mod$.
 
 b. If $\phi \in \od-mod$ is an open TCFT, then the closed state space of 
the open-closed field theory $L_\rho (\phi)$ (i.e., the value of the functor on the object given by the circle, $L_\rho (\phi) (S^1)$) is a chain
complex whose homology is given by the Hochschild homology of the $A_\infty$ -category, $\phi$.  That is,
$$
H_*(L_\rho (\phi)(S^1)) \cong HH_*(\cc_\phi).
$$
\end{theorem}

Here the Hochschild homology of a category enriched over chain complexes is computed via the Hochschild complex, whose
$n$-simplices are direct sums of terms of the form $Mor(\lambda_0, \lambda_1) \otimes Mor(\lambda_1, \lambda_2) \otimes \cdots \otimes Mor(\lambda_{n-1}, \lambda_n) \otimes  Mor(\lambda_n, \lambda_{0})$.   This is a double complex  whose boundary homomorphisms are the sum of the internal boundary maps in the chain complex of $n$-simplices,  plus the Hochschild boundary   homomorphism,  which is defined  as the alternating sum  $\sum_{i=0}^n (-1)^i  \p_i$,   where for $i=0, \cdots, n-1$,  $\p_i$ is induced by the composition
$Mor(\lambda_i, \lambda_{i+1}) \otimes Mor (\lambda_{i+1}, \lambda_{i+2}) \to Mor(\lambda_i, \lambda_{i+2})$.   $\p_n$ is induced by the composition
$Mor(\lambda_n, \lambda_{0}) \otimes Mor(\lambda_0, \lambda_{1}) \to Mor(\lambda_n, \lambda_{1})$.  
The Hochschild homology of an $A_\infty$-category enriched over chain complexes
is defined similarly.  See \cite{costello} for details.  

Costello's theorem can be interpreted as saying   that there is a ``universal" open-closed theory with a given value on the open cobordism category (i.e., the value of the derived left adjoint $L_\rho$), and that its closed state space has homology equal to the Hochschild homology of the associated $A_\infty$-category.  We note that in the interesting case when there is only one $\cd$-brane,  that is, $\cd = \{\lambda\}$, then
the $A_\infty$-category is an $A_\infty$-algebra,   and so the closed state space of the associated universal open-closed theory would have
homology given by the Hochschild homology of this algebra.    In particular this says that for any open-closed field theory $\phi$  with one   $\cd$-brane,  which has the corresponding $A_\infty$-algebra $A$,  then there is a well defined map from the Hochschild homology $HH_*(A) \to H_*(\phi (S^1))$.  This can be viewed as a derived version of the Moore-Segal result (Proposition \ref{center})  that gives a map $\phi (S^1) \to Z(A)$.  In the Moore-Segal setting,  $\phi (S^1)$ is an ungraded  Frobenius algebra ,  (or equivalently it has trivial grading) so we may identify it with $H_0(\phi (S^1))$.    Furthermore the center $Z(A)$ may be identified with the zero dimensional Hochschild \sl cohomology \rm $HH^0(A)$, so that the Moore-Segal result gives a map $H_0(\phi (S^1)) \to HH^0(A)$.  By the self-duality of the Frobenius algebra structures of $\phi (S^1)$ and of $A$, this gives a dual map $HH_0(A) \to H_0(\phi (S^1))$.  Costello's map can be viewed as  a derived version of this map.    

We end this section by remarking that recently Hopkins and Lurie have described a generalization of Costello's classification scheme that applies
in all dimensions.  The type of field theories they consider are called ``extended topological quantum field theories".  We refer the reader to \cite{lurie} for a description of their work.

\section{The string topology category and its Hochschild homology}

One of the goals of our project is to understand how string topology fits into Costello's picture.  The most basic operation in string topology  is the loop product defined by Chas and Sullivan \cite{chassullivan}:
$$
\mu : H_p(LM) \otimes H_q(LM) \to H_{p+q-n}(M)
$$
where $M$ is a closed, oriented, $n$-dimensional manifold.  Now let $B : H_q(LM) \to H_{q+1}(LM)$ be the operation induced by the rotation $S^1$-action on $LM$, $r : S^1 \times LM \to LM$. 
\begin{align}
B: H_q(LM)   &\to H_{q+1}(S^1 \times LM)  \xr{r_*} H_{q+1}(LM) \notag \\
\alpha &\to r_*([S^1]   \times   \alpha )   \notag
\end{align}

The following was one of the main theorems of \cite{chassullivan}.
\begin{theorem}\cite{chassullivan}  Let $\bh_*(LM) = H_{*+n}(LM)$ be the (regraded) homology of the free loop space.
Then with respect to the loop product $\mu$ and the degree one operator $B$,  $\bh_*(LM)$ has the structure of a (graded)
Batalin-Vilkovisky algebra.  That is, it is a graded commutative  algebra satisfying the following identities:
\begin{enumerate}
\item $ B^2 = 0,$   \quad \text{and}
\item
For $\alpha \in \bh_p(LM)$, and $\beta \in \bh_q(LM)$ the bracket operation 
$$
\{\alpha, \beta \} = (-1)^{|\alpha|}B(\alpha \cdot \beta) - (-1)^{|\alpha|}B(\alpha)\cdot \beta - \alpha \cdot B(\beta)
$$ is a derivation in each variable.  
 
\end{enumerate}
\end{theorem}
Moreover, a formal argument given in \cite{chassullivan} implies that the operation $ \{\,, \, \}$ satisfies   the (graded) Jacobi identity, and hence gives $\bh_*(LM)$ the structure of a graded Lie algebra.

The product is defined by considering the  mapping space, $Map (P, M)$ where $P$ is the pair of pants cobordism  (figure (\ref{figfour}))  between two circles and one circle.  By restricting maps to the incoming and outgoing boundaries, one has a correspondence diagram
\begin{equation}\label{inout}
LM  \xleftarrow{\rho_{out}} Map (P, M) \xr{\rho_{in}} LM \times LM.
\end{equation}
By retracting the surface $P$ to the homotopy equivalent   figure 8 graph, one sees that one has a homotopy cartesian square,
$$
\begin{CD}
Map (P, M)  @>\rho_{in} >>  LM \times LM \\
@VVV  @VVV \\
M @>>\Delta > M \times M 
\end{CD}
$$
where $\Delta : M \hk M \times M$ is the diagonal embedding.  This then allows the construction of an ``umkehr map"
$\rho_{in}^! : H_*(LM \times LM) \to H_{*-n}(Map (P,M))$.  This map was defined on the chain level in \cite{chassullivan},
and via a Pontryagin-Thom map $LM \times LM \to Map(P,M)^{TM}$ in \cite{cohenjones}.   Here $Map(P,M)^{TM}$ is the Thom
space of the tangent bundle $TM$ pulled back over the mapping space via evaluation at a basepoint, $Map(P,M) \to M$.   By twisting
with the virtual bundle $-TM$,  Cohen and Jones proved the following.
\begin{theorem}\cite{cohenjones}  For any closed manifold $M$, the Thom spectrum $LM^{-TM}$ is a ring spectrum.  When given an orientation of $M$,
the ring structure of $LM^{-TM}$ induces, via the Thom isomorphism, the Chas-Sullivan  algebra structure on $\bh_*(LM)$.
\end{theorem}

The Chas-Sullivan product was generalized to a TQFT by Cohen and Godin in \cite{cohengodin}.   Given a  cobordism  $\Sigma$ between
$p$-circles and $q$-circles,  they considered the following correspondence diagram analogous to (\ref{inout}).

\begin{equation}\label{correspond}
(LM)^q  \xleftarrow{\rho_{out}} Map (\Sigma, M) \xr{\rho_{in}}(LM)^p.
\end{equation}
Using fat (ribbon) graphs to model surfaces,  Cohen and Godin described an umkehr map
$$
\rho_{in}^! : H_*((LM)^p) \to H_{*+\chi (\Sigma)\cdot n}(Map (\Sigma, M))
$$
which allowed the definition of an operation 
$$
\mu_{\Sigma} = (\rho_{out})_* \circ \rho_{in}^! : H_*((LM)^p) \to  H_{*+\chi (\Sigma)\cdot n}((LM)^q)
$$
which yielded the (closed) TQFT structure.   In these formulae, $\chi (\Sigma)$ is the Euler characteristic of the cobordism $\Sigma$.  

Open-closed operations were first defined by Sullivan in \cite{sullivan}.  Somewhat later, Ramirez \cite{ramirez} and Harrelson \cite{harrelson}
showed that these operations define a positive boundary,  open-closed topological quantum field theory, in the Moore-Segal sense, except that the value of the theory lie in the category of graded vector spaces over a field $k$.    In this theory, which we call $\cs_M$,  the closed state space is given by
\begin{equation}\label{closedstate}
\cs_M (S^1) = H_*(LM; k).
\end{equation}
The set of $\cd$-branes $\cd_M$ is defined to be the set of connected, closed submanifolds $N \subset M$.  The value of this theory on
the interval labeled by submanifolds $N_1$ and $N_2$ (see figure (\ref{figfive}))  is given by
\begin{equation}\label{openpart}
\cs_M(N_1, N_2) = \cs_M(I_{N_1, N_2})   =  H_*(\cp_{N_1, N_2}),
\end{equation}
where $ \cp_{N_1, N_2 } $ is the space  of   paths $\alpha : [0,1] \to M$ with boundary conditions, $\alpha (0), \in N_1$, $\alpha (1) \in N_2$. 

Finally, using families of ribbon graphs modeling both closed and open-closed cobordisms,  in \cite{godin} Godin recently proved the following
result.

\begin{theorem}(Godin) \cite{godin}
Let $\co\cc_{\cd_M}^{H_*}$  be the category with the same objects as $\co \cc_{\cd_M}$, and whose morphisms are the homology of the morphisms in $\co \cc_{\cd_M}$.  That is,  given objects $c_1$ and $c_2$,  the morphisms from $c_1$ to $c_2$ are given by
$$
Mor_{\co\cc_{\cd_M}^{H_*}}(c_1, c_2) = H_*(\cm_{\cd} (c_1, c_2);k)  \cong   \bigoplus H_*(BDiff^+ (\Sigma, \p \Sigma);k )
$$
where the direct sum is taken over all diffeomorphism classes of open-closed cobordisms from $c_1$ to $c_2$.   
Then the above string topology operations can be extended to a symmetric monoidal functor
$$
\cs_M : \co\cc_{\cd_M}^{H_*}  \to Gr \, Vect
$$
where $Gr \, Vect$ is the category of graded vector spaces over $k$, whose monoidal structure is given by (graded) tensor product.
In other words, the string topology of $M$ is a  positive boundary, open-closed  ``homological conformal field theory" (HCFT).
\end{theorem}

Notice that being a homological conformal field theory is a weaker property than being a topological conformal field theory, and so Costello's results cannot be immediately applied to the string topology of a manifold $M$.  In order for the functor $\cs_M$ to actually induce a TCFT, the string topology operations must be defined on the chain level, and satisfy the appropriate compatibility and coherence properties. It is conjectured
that in fact this can be done.  In any case, Costello's theorem (Theorem \ref{costello} above) suggests  that there is an $A_\infty$-category associated to the string topology of $M$, and that its Hochschild homology should be the value of the closed state space, $\cs_M(S^1) = H_*(LM; k)$.   Theorem \ref{bct1} in the introduction asserted the existence of such a category; we will describe the construction in more detail below, although full proofs appear in \cite{bct}.

Another interesting question arises when there is only a single
$D$-brane $\cd = \{N\}$, where $N$ is a fixed, connected submanifold
of $M$.  In this case the corresponding $A_\infty$-category would be
an $A_\infty$-algebra.  Here it turns out that for Poincar\'{e} duality
reasons it is more appropriate to consider Hochschild cohomology.  The
question described in the introduction, regarding the relationship
between these Hochschild cohomology algebras and the Chas-Sullivan
algebra structure on $H_*(LM)$, was based on the idea that string
topology, even in this ``one $D$-brane" setting should fit into
Costello's picture of a universal open-closed TCFT.  In particular the
calculations described below verify that for a large class of
submanifolds $N \subset M$, the full subcategory of $\cs_M$ consisting
of the single object $N$ still yields the full closed state space of
string topology,
$$ HH^*(C_*(\cp_{N,N}), C_*(\cp_{N,N})) \cong H_*(LM) = \cs_M(S^1).
$$
 
  An important idea that runs throughout the proofs of these
  statements is that of a ``derived" form of Poincar\'{e} duality.
  Namely, instead of the classical setting where one has coefficients
  given by modules over the group ring $\bz [\pi_1(M)]$, we need  a
  version of Poincar\'{e} duality that applies to modules over the
  differential graded algebra $C_*(\Omega M)$.  To be precise, what we
  mean by $C_*(\Omega M)$ is the DGA corresponding to the $Hk$-module
  spectrum $Hk \sma \Sigma^{\infty}  (\Omega M_+)$ via the equivalence
  of \cite{shipleyHZ}, where here $\Omega M$ denotes a model of the
  based loop space which is a topological group, and $Hk$ is the Eilenberg-MacLane spectrum for the field $k$. 
 
 More specifically,  recall that for any connected space $X$, one has
 natural Eilenberg-Moore equivalences, 
 \begin{align}
 Tor_{C_*(\Omega X)} (k, k) &\cong H_*(X; k) \notag \\
 Ext_{C_*(\Omega X)} (k, k) &\cong H^*(X; k). \notag
 \end{align}
 (Here and below we will suppress grading.)  Written on the level of
 chain complexes there are equivalences, 
 \begin{align}
k \otimes^L_{C_*(\Omega X)} k &\simeq C_*(X; k) \notag \\
 Rhom_{C_*(\Omega X)} (k, k) &\simeq C^*(X; k) \notag
 \end{align}
 the later equivalence being one of DGA's.  
 Now suppose $M$ is a connected, closed, oriented manifold (or more generally,  a Poincar\'{e} duality space).  Then the fundamental class,
 $[M] \in  Tor_{C_*(\Omega M)} (k, k) \cong H_*(M; k)$ is represented by a fundamental cycle $[M] \in k \otimes^L_{C_*(\Omega M)} k$.  Classical Poincar\'{e} duality can be viewed as saying that capping induces a chain homotopy equivalence,
 $$
 \cap [M] :  Rhom_{C_*(\Omega M)} (k, k)  \xr{\simeq} k \otimes^L_{C_*(\Omega M)} k.
 $$
 It follows from the work of Dwyer-Greenlees-Iyengar \cite{DGI} and Klein \cite{klein}  that indeed one has a chain equivalence,
 \begin{equation}\label{PD}
  \cap [M] :  Rhom_{C_*(\Omega M)} (k, P)  \xr{\simeq} k \otimes^L_{C_*(\Omega M)} P,
  \end{equation}
  where $P$ is \sl any \rm differential graded module over $C_*(\Omega M)$, that is bounded below. 
      This is the ``derived" form of Poincar\'{e} duality that we need.
      
\begin{remark}      
This interpretation of the results of \cite{DGI} and \cite{klein} was used by  Malm in \cite{malm} to study other aspects of string topology in this algebraic setting.
\end{remark}

In what follows we sketch how this duality is used to construct the
string topology category $\cs_M$.  Given a connected, closed, oriented
submanifold $N \subset M$, let $F_N$ be the homotopy fiber.  Because of our connectivity hypotheses, one can view $F_N$ as
$\cp_{N, x_0}$, where $x_0 \in M$ is a fixed basepoint.  We then have
homotopy fibrations, 
\begin{align}
F_N \to N \to M \notag \\
\Omega M \to F_N \to N \notag \\
\Omega N \to \Omega M \to F_N  \notag 
\end{align}
(In the last sequence, we   choose a specific equivalence between the homotopy fiber of $\Omega M \to F_N$ and $\Omega N$. The choice is not canonical since we are not assuming that the inclusion $N \hk M$ preserves basepoints.    This choice, however,  is not used in our definition of the string topology category below.)

From these fibrations we can regard $C_*(F_N)$ as a $DG$-module over
$C_*(\Omega M)$.  Once again, what we mean by $C_*(F_N)$ is the
$DG$-module corresponding to the module spectrum $Hk \sma
\Sigma^{\infty} ((F_N)_+)$.

A key observation is that the compactness of $N$ and $M$ impose a
strong condition on $C_*(F_N)$.  Specifically, we have the following
lemma.

\begin{lemma}
Let $N$ and $M$ be finite-dimensional connected  complexes with $N
\subset M$,  and $F$ the homotopy fiber of the inclusion $N \to M$.
Then $C_*(F)$ is small as a $C_* (\Omega M)$-module.
\end{lemma}

There are equivalences
\begin{align}\label{fnequiv}
k \otimes^L_{C_*(\Omega M)} C_*(F_N) \simeq C_*(N) \notag \\
 k  \otimes^L_{C_*(\Omega N)} C_*(\Omega M)  \simeq C_*(F_N).
\end{align}
(The second of these equivalences depends on the choice of equivalence of $\Omega N$ and the homotopy fiber of $\Omega M \to F_N$ above.)
Furthermore, we have homotopy cartesian squares
$$
\begin{CD}
\cp_{N_1, N_2}  @>>>  N_2  \\
@VVV  @VVV \\
N_1  @>>>  M
\end{CD}
\quad \text{and}  \quad
\begin{CD}
\Omega M  @>>> F_{N_2}  \\
@VVV  @VVV \\
F_{N_1}  @>>>  \cp_{N_1, N_2}.
\end{CD}
$$
Therefore an {\em Eilenberg-Moore} argument yields a chain homotopy equivalence
\[
C_*(\cp_{N_1, N_2}) \simeq C_*(F_{N_1}) \otimes^L_{C_*(\Omega M)}
 C_*(F_{N_2}).
\]
We can then make the following equivalences:
  \begin{align}
  C_*(\cp_{N_1, N_2}) &\simeq  
 C_*(F_{N_1}) \otimes^L_{C_*(\Omega M)}   C_*(F_{N_2}) \notag \\
  & \simeq k \otimes^L_{C_*(\Omega N_1)} C_*(F_{N_2})  \quad \text{by (\ref{fnequiv}) and change of rings}, \notag \\
 &\simeq Rhom_{C_*(\Omega N_1)} (k, C_*(F_{N_2})) \quad \text{by
   Poincar\'{e} duality equivalence (\ref{PD}) for } N_1,\notag \\
 &\simeq Rhom_{C_*(\Omega M)} (C_*(F_{N_1})  , C_*(F_{N_2})) \quad \text{again by (\ref{fnequiv}) and change of rings}. \notag
 \end{align}
 
 We remind the reader that in the above equivalences all gradings (and
 grading shifts) are suppressed.  Further, it is worth emphasizing
 that it is Poincar\'{e} duality for $N_1$ that is used in these
 equivalences.
 
By using cofibrant-fibrant replacements of $C_* (F_{N_i})$ (for which we
use the standard model structure on $DG$-modules over $C_*(\Omega M)$,
e.g., \cite[\S 7]{johnson}, \cite{schwede-shipley, shipleyHZ}), we can
regard the derived homomorphism complexes as possessing a strict
composition pairing.  This observation gives rise to the definition of
the string topology category.
 
\begin{definition}\label{definition}  Let $M$ be a connected, closed, oriented manifold with fixed basepoint $x_0 \in M$.  
The string topology category $\cs_M$ has as
\hspace{5 pt}
\begin{enumerate} 
\item Objects the pairs $(N, F_N)$, where $N$ is a  connected,  closed, oriented submanifold $N \subset M$ and $F_N$ is a specific choice of model for the homotopy fiber of $N \to M$ with an action of $\Omega M$.
\item Morphisms from $N_1$ to $N_2$ the derived homomorphism complex
\[Rhom_{C_*(\Omega M)}(C_*(F_{N_1}), C_*(F_{N_2})),\] computed via
functorial cofibrant-fibrant replacement of $C_*(F_{N_i})$. 
\end{enumerate}
In other words, $\cs_M$ is the full subcategory of the $DG$-category
of differential graded modules over $C_*(\Omega M)$ with objects
cofibrant-fibrant replacements of $C_*(F_{N})$ for $N \subset M$ a
submanifold as above.
\end{definition}

We remark that this derived form of Poincar\'{e} duality has another
interpretation, namely one of a derived and non-commutative analogue
of local Serre duality for a regular, $n$-dimensional local ring
$R$. For any $R$-module $S$, we have the isomorphism
\begin{equation}\label{iso}
\mathrm{Rhom}_R(k,S) \simeq k\otimes_R S[-n].  
\end{equation}
Indeed, the Gorenstein property says that 
\[
\mathrm{Ext}^n_R(k,R) = k,\qquad  \mathrm{while\ the\ other\ Ext's\ are\ zero,} 
\]
and (\ref{iso}) follows by standard homological algebra arguments and 
finiteness of $k$ over $R$ (a consequence of regularity).

The Gorenstein property of the DGA $R = C_*(\Omega{N})$ follows from 
Poincar\'{e} duality for $N$ and from the Eilenberg-Moore isomorphism 
\[
\mathrm{Rhom}_{C_*(\Omega{N})}(k,C_*(\Omega{N})) \simeq 
C^*(N; \tilde{C}_*(\Omega{N})),
\]
where $\tilde{C}_*(\Omega{N})$ is the cohomology-coefficient system over 
$N$ defined by the fibrewise chains of the path fibration $PN \to N$. Indeed, 
Poincar\'{e} duality identifies the last cohomology with the total homology 
of the based path space, shifted up by $n$, and we conclude that the right-hand 
complex is isomorphic to $k$ placed in degree $n$.  Regularity of
$C_*(\Omega{N})$ (i.e., finiteness of $k$ as a $C_*(\Omega N)$-module)
is a consequence of the fact that $N$ has a finite cell decomposition.  See
\cite[\S 10]{DGI} for further discussion of the relationship of the
Gorenstein condition to topological duality phenomena. 
  
Finally, we point out that this generalized notion of Poincar\'{e}
duality (i.e., with coefficients being modules over $C_*(\Omega M)$)
is at the heart of the argument that shows that composition in this
category realizes the string topology compositions on the level of
homology,
$$ 
H_*(\cp_{N_2, N_3}) \otimes H_*(\cp_{N_1, N_2}) \to H_*(\cp_{N_1, N_3}).
$$

Curiously, this is most easily seen on the level of spectra.  As will
be shown in \cite{bct}, if one runs through the above chain level
argument on the level of module spectra over ring spectra, then one
proves a twisted form of Atiyah duality: 
$$
Rhom_{\Sigma^\infty (\Omega M_+)}(\Sigma^\infty (F_{N_1}),
\Sigma^\infty (F_{N_2})) \simeq \cp_{N_1, N_2}^{-TN_1} 
$$
where $\cp_{N_1, N_2}^{-TN_1}$ is the Thom spectrum of the virtual bundle $-TN_1$ obtained by pulling back the negative tangent bundle $-TN_1 \to N_1$
over $\cp_{N_1, N_2}$ via the map $ \cp_{N_1, N_2} \to N_1$ that evaluates a path at its starting point.  We view this as a twisted form of Atiyah duality, because in the case when $N_1 = N_2 = M$, then $F_M \simeq point$, and $\cp_{M,M} \simeq M$.  We then  have the traditional form of Atiyah duality:
$$
Rhom_{\Sigma^\infty (\Omega M_+)} (S^0, S^0)   \simeq  M^{-TM}.
$$
Here $S^0$ is the sphere spectrum, and $Rhom_{\Sigma^\infty (\Omega
  M_+)} (S^0, S^0) \simeq Rhom_{S^0}(\Sigma^\infty (M_+), S^0)$ is the
Spanier-Whitehead dual of $M$.   This kind of twisted Atiyah duality
has been studied before in several contexts.  For example, in
\cite{umkehr},  it was studied in the context of ``Poincar\'{e} duality
with coefficients in a fibered spectrum".  It also appeared, in a
slightly different context in \cite{klein}, and in \cite{Hu}.     In
the classical setting, as well as these ``twisted settings", it is
fairly standard to see that the composition pairings correspond, up to
homotopy,  to pairings constructed on the level of  Thom spectra via
the Pontrjagin-Thom construction (see \cite{atiyahdual} for this type
of argument).  However the Pontrjagin-Thom constructions are precisely
how the string topology pairings are defined, as in \cite{cohenjones}
and \cite{ramirez}.   When one passes to chains, this implies that the
composition of morphisms in the string topology category, correspond,
on the level of homology, to the open-closed string topology
operations.   Details of these arguments will appear in \cite{bct}.

When $M$ is simply-connected, we have a useful ``dual'' model of
$\cs_M$ which follows from the following alternative description of
the path spaces $C_*(\cp_{N_0, N_1})$, another consequence of the
generalized form of Poincar\'{e} duality.  Here we regard $C^*(M)$ as an
$E_\infty$-algebra and we are relying on the existence of a model
structures on modules over an $E_\infty$-algebra \cite{mandell}.

\begin{lemma}  Let $M$ be a closed, simply connected manifold, and
  $N_0, N_1 \subset M$ connected, oriented, closed submanifolds.  Then
  there is a chain equivalence, 
 $$ C_*(\cp_{N_0, N_1}) \simeq Rhom_{C^*(M)} (C^*(N_1), C^*(N_0)).
 $$ Equivalently, there is a chain equivalence
 $$ Rhom_{C_*(\Omega M)}(C_*(F_{N_0}), C_*(F_{N_1})) \simeq
  Rhom_{C^*(M)} (C^*(N_1), C^*(N_0)).
 $$
 \end{lemma}
 
Furthermore, it turns out that the multiplicative structures are
compatible with the equivalences of the previous lemma and as a
consequence, we obtain the following comparison result.  Recall that
when comparing enriched categories, the correct notion of equivalence
is given by considering enriched functors which induce equivalences on
mapping objects and an underlying equivalence of ``homotopy
categories'' (e.g., defined by passing to $H^0$ for $DG$-categories or
components for spectral categories).  We will refer to such
equivalences in general as Dwyer-Kan equivalences, although in special
cases they tend to have specific names (e.g., quasi-equivalence of
$DG$-categories).
  
\begin{theorem} For $M$ a closed, simply connected manifold, there is
  a zigzag of Dwyer-Kan equivalences between the string topology category
  $\cs_M$ and the full subcategory of the category of $C^*(M)$-
  submodules, with objects cofibrant-fibrant replacements (as
  $C^*(M)$-modules) of the cochains $C^*(N)$ for $N \subset M$ a
  connected, oriented submanifold.
\end{theorem}

The Hochschild homology statement in Theorem~\ref{bct1} follows from
an identification $HH_*(\cs_M) \cong HH_*(\Omega M)$.  This in turn is
obtained as a straightforward consequence of the general theory
developed in \cite{blumbergmandell}; the thick closure of $\cs_M$
inside the category of $C_*(\Omega M)$-modules is the entire category
of finite $C_*(\Omega M)$-modules.  This is essentially a Morita
equivalence result. When one restricts to a single object  $N$, (the ``one-brane" situation),   the endomorphism algebra
is equivalent to $C_*(\cp_{N,N})$.  In this case analysis of the Hochschild (co)homology  requires a more
involved  Morita theory.  This is because it is definitely not 
the case in general
that there is a Morita equivalence between the category of $C_*(\Omega
M)$-modules and the category of $C_*(\cp_{N,N})$-modules.

Here the situation is a kind of Koszul duality, and so whereas the
Hochschild homologies of the path algebras vary as $N$ varies
(explicit descriptions will be given in \cite{bct}), it is reasonable to expect that the
Hochschild cohomologies should coincide. 
To approach these calculations,  the following basic principle is used in \cite{bct}:

\begin{theorem}  Let  $R$ and $S$   be  two differential graded algebras over a field $k$,  and suppose there exist  $R-S$  (differential graded) modules satisfying the following equivalences:
\begin{equation}\label{morita}
 Rhom_{R}(P, Q) \simeq S  \quad \text{and} \quad Rhom_{S}(P, Q) \simeq R.
\end{equation}
 Then their Hochschild   cohomologies are isomorphic,
 $$
 HH^*(R, R) \cong HH^*(S, S).
 $$
\end{theorem}
 
Using this result, given a submanifold $N \subset M$,  one considers $$R \simeq C_*(\Omega M),  \quad   S \simeq C_*(\cp_{N, N}) \simeq Rhom_{C_*(\Omega M)}(C_*(F_N), C_*(F_N)).$$   The $R-S$ modules are both given by $P = Q \simeq C_*(F_N)$.     We already know that
 $ Rhom_{R}(P, Q)    \simeq S$, for any $N \subset M$,  when $M$ is simply connected.  The Hochschild cohomology calculations are then reduced 
to a question about ``double centralizers'':  
$HH^*(C_*(\cp_{N,N}))$ is equivalent to $HH^*(C_*(\Omega M))$ if there
is an equivalence
\[Rhom_{C_*(\cp_{N,N})}(C_*(F_N), \,  C_*(F_N)) \simeq C_*(\Omega
M).\]
Keller \cite{keller} has shown a more general formulation of the
sufficiency of a double-centralizer condition for the equivalence of
Hochschild cohomology for $DG$-categories (and note that there is also
a generalization of his theorem to $THH$ cohomology and spectral
categories using the technology of \cite{blumbergmandell}).

The question of the existence of such equivalences can be studied
using the generalized Morita theory of Dwyer, Greenlees, and Iyengar
\cite{DGI}, and leads to the following characterization:
 
\begin{theorem}\label{doublecent}  Assume $M$ is simply connected.  Then for any $N \subset M$ in $\cd$, 
 $$ Rhom_{C_*(\cp_{N,N})}(C_*(F_N), \, C_*(F_N)) \simeq \hat
  C_*(\Omega M)
 $$ Here $\hat C_*(\Omega M)$ is the Bousfield localization of
  $C_*(\Omega M)$ with respect to the homology theory $h_*^N$, defined
  on the category of $C_*(\Omega M)$-modules given by
 $$ h_*^N(P) = Ext_{C_*(\Omega M)}(C_* (F_N), P) \cong Ext_{C_*(\Omega N)}(k, P).
 $$
(Note that this is best regarded as a completion process, despite the
  terminology of localization; we will refer to local objects as
  $C_*(F_N)$-complete.)
\end{theorem}

An immediate corollary is that the double centralizer property holds
if and only if $C_* (\Omega M)$ is $C_* (F_N)$-complete.  Therefore,
we are immediately led to study the following question: For which
submanifolds $N \subset M$ is $C_*(\Omega M)$ in fact
$C_*(F_N)$-complete?  Counterexamples (obtained in consultation with
Bill Dwyer) exist that suggest that this does not always hold.
However, we can show the result in certain useful special cases.  In
particular, we know that $C_*(\Omega M)$ is $C_*(F_N)$-complete in the
following cases:

 \begin{enumerate}
 \item The inclusion map $N \hk M$ is null homotopic. This implies
   that $F_N \simeq \Omega M \times N$, and $\cp_{N,N} \simeq \Omega M
   \times N \times N$.
 \item The inclusion $N \hk M$ is the inclusion of the fiber of a
   fibration $p : M \to B$.  More generally there is a sequence of
   inclusions, $N \subset N_1 \subset N_2 \subset \cdots \subset N_k =
   M$ where each $N_i \subset N_{i+1}$ is the inclusion of the fiber
   of a fibration $p_{i+1} : N_{i+1} \to B_{i+1}. $
 \end{enumerate}
 
These results also have consequences for certain module categories
related to the rings we are considering.  Denote by $E_N$ the
endomorphism ring $Rhom_{C_*(\Omega M)}(C_*(F_{N}), C_*(F_{N}))$,
which as noted above provides a strictly multiplicative model of
$C_*(\cp_{N,N})$.  

\begin{theorem}\label{module} 
When the double centralizer condition holds, the following categories
of modules are Dwyer-Kan equivalent:
\begin{enumerate}
\item The thick subcategory of $C^*(M)$-modules generated by $C^*(M)$
  (i.e., the perfect modules).
\item The thick subcategory of $C_*(\Omega M)$-modules generated by
  the trivial module $k$.
\item The thick subcategory of $E_N$-modules generated by
  $Rhom_{C_*(\Omega M)}(C_*(F_N), k)$.  Notice that this latter module is
  equivalent to $C^*(N)$.
\end{enumerate}
\end{theorem}

Note that the equivalence of (1) and (2) was shown in \cite{DGI};
their methods extend to provide the comparison to (3).     These equivalences of categories of modules
are relevant to the Floer theory of compact Lagrangians in  the cotangent bundle $T^*M$, as we
discuss in Section 3.
 
Finally, we point out that there are spectrum level analogues of the
above theorems (with essentially similar proofs), in particular
Theorem  \ref{bct1} as stated in the introduction.
  
\begin{theorem}\label{spectrum}  There is a string topology category,
   enriched over spectra, which by abuse of notation we still refer to
   as $\cs_M$, whose objects again are determined by elements of $\cd_M$.  The
   morphism spectrum between $N_1$ and $N_2$ is the analogous mapping
   spectrum
\[Rhom_{\Sigma^{\infty}_+ \Omega M}(\Sigma^{\infty}_+ F_{N_1},
\Sigma^{\infty}_+ F_{N_2}),\]
 and has the homotopy type of  
 $\cp_{N_1, N_2}^{-TN_1}$, the Thom spectrum of the virtual bundle
   $-TN_1$, where $TN_1$ is the tangent bundle of $N_1$, pulled back
   over $\cp_{N_1,N_2}$ via the evaluation map that takes a path
   $\alpha \in \cp_{N_1,N_2}$ to its initial point $\alpha (0) \in
   N_1$.  Furthermore, the $THH$ of this category is
   the suspension spectrum of the free loop space (with a disjoint
   basepoint),
  \begin{equation}\label{thhcat}
 THH_\bullet (\cs_M) \simeq \Sigma^\infty (LM_+).
 \end{equation}
 \end{theorem}  
 
 Moreover, the analogue of the above Hochschild cohomology statement
 is the following:
 \begin{theorem}\label{onebranespec}
Assume that $M$ is a simply connected, closed manifold.  Then given
any connected, closed, oriented submanifold $N \subset M$ for which
the double-centralizer condition holds, then the Thom spectrum
$\cp_{N,N}^{-TN}$ is a ring spectrum, and its topological Hochschild
cohomology is given by
 \begin{equation}\label{onebranesp}
 THH^\bullet (\cp_{N,N}^{-TN}) \simeq LM^{-TM}
\end{equation}
 and the equivalence is one of ring spectra.
 \end{theorem}

\section{Relations with the Fukaya category of the cotangent bundle}   
 
This section is speculative, regarding the possible relationships between the string topology category $\cs_M$, and the Fukaya category of the cotangent bundle          $T^*M$.  The Fukaya category is an $A_\infty$- category associated to a symplectic manifold $(N^{2n}, \omega)$. Here $\omega \in \Omega^2(N)$ is a symplectic $2$-form.  Recall that for any smooth $n$-manifold $M^n$, $T^*M$ has the structure of an \sl exact \rm symplectic manifold.  That is, it has a symplectic $2$-form $\omega$ which is exact.  In the case $T^*M$,  $\omega = d\theta$, where $\theta$ is the Liouville one-form defined as follows.  Let $p : T^*M \to M$ be the projection map.  Let $x \in M$, and $t \in T^*_xM$.  Then 
  $\theta (x,t)$ is the given by the composition,
 $$
 \theta (x,t): T_{x,t}(T^*M) \xr{dp}T_xM \xr{t}\br.
 $$
 There has been a considerable amount of work comparing the symplectic topology of $T^*M$ with the string topology of $M$.  This relationship begins with 
 a theorem of Viterbo \cite{viterbo}, that the symplectic Floer homology is isomorphic to the homology of the free loop space,
 $$
 SH_*(T^*M) \cong H_*(LM).
 $$
 The symplectic Floer homology is computed via a Morse-type complex   associated to the (possibly perturbed) ``symplectic action functional",  $\ca : L(T^*M) \to \br$.
 The perturbation is via a choice of Hamiltonian, and so long as the Hamiltonian grows at least quadratically near infinity, the symplectic Floer homology
 is described by the above isomorphism.   The precise relationship between the Floer theory of the symplectic action functional $\ca$ and Morse theory on $LM$ was studied in great detail by Abbondandolo and Schwarz in \cite{abschwarz}.  In particular they were able to show that a ``pair of pants" (or ``quantum") product construction in $SH_*(T^*M)$ corresponds under this isomorphism to a Morse-theoretic analogue of the Chas-Sullivan product in $H_*(LM)$.  In \cite{cohenschwarz} this
 product was shown to agree with the Chas-Sullivan construction.

 The objects of the Fukaya category $Fuk(T^*M)$ are exact, Lagrangian submanifolds $L \subset T^*M$.  The morphisms are the ``Lagrangian intersection Floer cochains",
 $CF^*(L_0, L_1)$.  These Floer cochain groups are also a Morse type cochain complex  associated to a functional on the path space, 
 $$\ca_{L_0, L_1} : \cp_{L_0, L_1}(T^*M) \to \br. $$
 If $L_0$ and $L_1$ intersect transversally, then the critical points are the intersection points (viewed as constant paths), and the coboundary homomorphisms
 are computed by counting holomorphic disks with prescribed boundary conditions.  Of course if $\ca_{L_0, L_1}$ were actually a Morse function, satisfying the Palais-Smale convergence conditions, then these complexes would compute $H_*(\cp_{L_0, L_1}(T^*M))$.  One knows that this Morse condition is \sl not \rm satisfied,  but there are examples, when this homological consequence is nonetheless satisfied.  Namely,  let $N \subset M$ be an oriented, closed submanifold.  Let $\nu_N$ be the conormal bundle.  That is, for $x \in N$,  $\nu_N (x) \subset T^*_x M$ consists of those cotangent vectors which vanish on the subspace $T_xN \subset T_xM$.  Notice that the conormal bundle is always an $n$-dimensional submanifold of the $2n$-dimensional manifold $T^*M$.  It is a standard fact that  the conormal bundle $\nu_x N$ is a (noncompact) Lagrangian submanifold of $T^*M$.  Notice that for any two closed, oriented submanifolds   $N_0, N_1 \subset M$ the following  path spaces in the cotangent bundle and  in the base manifold $M$ are homotopy equivalent:
 $$
 \cp_{\nu_{N_0},\nu_{N_1}}(T^*M))  \simeq \cp_{N_0, N_1} (M).
 $$
 The following was recently proven by Abbondandalo, Portaluri, and Schwarz \cite{abboschwarz3}:
 
 \begin{theorem}   Given any closed, oriented submanifolds $N_0, N_1 \subset M$, then the intersection Floer cohomology of the $HF^*(\nu_{N_0}, \nu_{N_1})$ is isomorphic to the homology of the path space,
 $$
 HF^*(\nu_{N_0}, \nu_{N_1})  \cong H_*(\cp_{N_0, N_1}(M)).
 $$
 \end{theorem}
 
 If one could realize these isomorphisms on the level of chain complexes, in such a way that the compositions correspond,  then one would have a proof of the following conjecture:
 
 \begin{conjecture}   Let $Fuk_{conor}(T^*M)$ be the full subcategory of the Fukaya category generated by conormal bundles of closed, connected, submanifolds of $M$.  Then there is a Dwyer-Kan equivalence with the string topology category,
 $$
 Fuk_{conor}(T^*M) \simeq \cs_M.
 $$
 \end{conjecture}
 
\begin{remark}
\hspace{5 pt}
 \begin{enumerate}
 \item In this conjecture one probably wants to study the ``wrapped" Fukaya category as defined by Fukaya, Seidel, and Smith in \cite{fuk-seid-smith}.   
 \item  When $N = point$, its conormal bundle is a cotangent fiber, $T_x M$.   Abouzaid \cite{abouzaid} has recently described the $A_\infty$ relationship between the Floer
 cochains $CF^*(T_x M, T_x M)$ and the chains of the based loop space, $C_*(\Omega M)$.  
\end{enumerate} 
\end{remark}
 
 There are other potential relationships between the Fukaya category
 and the string topology category as well. 
 For example, Fukaya, Seidel, and Smith \cite{fuk-seid-smith}  as well as Nadler \cite{nadler} building on work of Nadler and Zaslow \cite{nadlerzas} showed that when $M$ is simply connected, the Fukaya category $Fuk_{cpt}(T^*M)$ generated by compact, exact Lagrangians with Maslov index zero, has a fully faithful embedding into the derived category of modules over the Floer cochains, $CF^*(M, M)$   where $M$ is viewed as a Lagrangian submanifold of $T^*M$ as the zero section.  Furthermore, one knows that the Floer cohomology, $HF^*(M, M)$ is isomorphic to $H^*(M)$, and recently
 Abouzaid \cite{abouzaid2} proved that $CF^*(M,M) \simeq C^*(M)$ as $A_\infty$- differential graded algebras.  So $Fuk_{cpt}(T^*M)$ can be viewed as a sub-$A_\infty$-category of the derived category of $C^*(M)$-modules.  
 
When $M$ is simply connected, recall that the string topology category
$\cs_M$ can also be viewed as a subcategory of the category of
$C^*(M)$- modules.  From Nadler's work we see that the relationship
between the compact Fukaya category $Fuk_{cpt}(T^*M)$ should be
equivalent to the ``one-brane" string topology category, $\cs_M^M$,
which is the full subcategory  of $\cs_M$ where the only $D$-brane
is the entire manifold itself.  This category, in turn, generates the
derived category of perfect $C^*M$-modules.
 
 The significance of this potential relationship is amplified when one
 considers recent work of Hopkins and Lurie \cite{lurie} classifying
 ``extended" topological conformal field theories.  This can be viewed
 as a direct generalization of the work of Moore-Segal, and of
 Costello discussed above.  In their classification scheme, such a
 field theory is determined by an appropriately defined ``Calabi-Yau"
 category enriched over chain complexes.  The category of perfect
 $C^*(M)$-modules is such a category.  Moreover generalized Morita theory
 implies that this category is Dwyer-Kan equivalent to the category of
 $k$-finite $C_*(\Omega M)$-modules.  Thus these categories should
 determine an extended field theory, which should correspond to string
 topology.  On the other hand, by the above remarks, the Fukaya
 category $Fuk_{cat}(T^*M)$ should also determine a field theory,
 presumably the ``Symplectic Field Theory" of Eliashberg, Givental,
 and Hofer \cite{EHG} applied to $T^*M$.  One can therefore speculate
 that this line of reasoning may produce an equivalence of the
 symplectic field theory of $T^*M$, and of the string topology of $M$.
 There is evidence that such an equivalence may exist, for example the
 work of Cielebak and Latchev \cite{cielelai}.  Pursuing this
 relationship using the Hopkins-Lurie classification scheme could lead
 to a very satisfying understanding of the deep connections between
 these two important theories.

 \end{document}